\documentclass[smallextended,nospthms]{svjour3}  
\usepackage{latexsym} 
\usepackage{amssymb}
 \usepackage{amsthm}
 \usepackage{amsmath}
 \numberwithin{equation}{section}

\def\bbb#1{\hbox {{\gordas #1}}}

\font\gordas = msbm10 at 12pt
\def\bbb#1{\hbox {{\gordas #1}}}

\def\ene{{\bbb N}}
\def\UNO{1\mkern-7mu1}

\def\ene{{\bbb N}}

\theoremstyle{plain}
\newtheorem{thm}{Theorem}[section]
\newtheorem{prop}[thm]{Proposition}
\newtheorem{lem}[thm]{Lemma}
\newtheorem{cor}[thm]{Corollary}
\theoremstyle{definition}
\newtheorem{defn}[thm]{Definition}
\newtheorem{rem}[thm]{Remark}

\newcommand{\Pn}{\mathcal{P}(n)}
\newcommand{\G}{{\mathcal {G}}}
\newcommand{\Ss}{{\mathcal {S}}}
\newcommand{\R}{\mathbb{R}}
\newcommand{\abs}[1]{\left|#1\right|}

\begin{document}

\title{ From intersection local time to the Rosenblatt process
\thanks{Partially supported by Conacyt grant 98998 (Mexico) and NCN grant 
2012/07/B/ST1/03417 (Poland).
}
}


\author{T. Bojdecki \and Luis G. Gorostiza \and Anna Talarczyk
}


\institute{T. Bojdecki \at
University of Warsaw, Warsaw, Poland\\
e-mail: tobojd@mimuw.edu.pl
\and
L.G. Gorostiza \at
Centro de Investigaci\'on y de Estudios Avanzados, Mexico City, 
Mexico\\
e-mail: lgorosti@math.cinvestav.mx
\and
A. Talarczyk \at
University of Warsaw, Warsaw, Poland\\
e-mail: annatal@mimuw.edu.pl}

\date{Received: date / Accepted: date}

\maketitle

\begin{abstract}
The Rosenblatt process was obtained by Taqqu (1975) from convergence in 
distribution of partial sums of strongly dependent random variables. In 
this paper
we give a particle picture approach to the Rosenblatt  process with the 
help of 
intersection  local time and white noise analysis, and discuss
measuring its long range dependence by means of a number called 
dependence exponent.
\keywords{Rosenblatt process \and  long range dependence \and intersection 
local time \and particle system.}
\subclass{MSC Primary 60G18 \and Secondary 60F17}
\end{abstract}
%

\section{Introduction}
\setcounter{section}{1}
\setcounter{equation}{0}

The Rosenblatt process was introduced and studied by Taqqu \cite{T1}, 
motivated by a counter-example of Rosenblatt regarding a strong mixing 
condition \cite{R}. See \cite{T4} for the history and an overview of 
the 
process, and related work. The Rosenblatt process arises from convergence 
in distribution of normalized partial sums of strongly dependent random 
variables. It lives in 
the second Wiener chaos, in contrast to Gaussian processes (which belong 
to the first chaos) that are obtained from sums of independent or weakly 
dependent random variables. The Rosenblatt process possesses some of the 
main properties of the (Gaussian) fractional Brownian motion with Hurst 
parameter greater than $1/2$: mean zero, H\"older continuity, non differentiability, selfsimilarity
 and stationarity of increments (hence it has the same form of covariance as 
 fractional Brownian motion), infinite divisibility, long range dependence in the sense that 
the sum over $k$ of covariances of increments on the 
 intervals $[0,1]$ and 
$[k,k+1]$ diverges, and it is not a semimartingale. The Rosenblatt process 
is the simplest non Gaussian Hermite process (see \cite{T2}).
 It is a counterpart of the 
fractional Brownian motion, which is the most prominent long range 
dependent Gaussian process. Significant attention has been attracted by 
the Rosenblatt process due to its mathematical interest, and possible 
applications where the Gaussian property may not be assumed. Recent papers 
on the subject include \cite{Tu}, which develops a related stochastic 
calculus  and mentions areas of application (see also references therein), 
\cite{MT1} and \cite{VT}, where new properties of the process have been 
found, and \cite{GTT}, which provides a strong approximation for the process. 
Relevant information for the present paper on the Rosenblatt process and fractional Brownian 
motion is given in the next section.

Our main objective in this paper is to show a different way of obtaining 
the Rosenblatt process. The method of Taqqu, which was developed for 
Hermite processes generally \cite{T2}, is based on limits of 
sums of strongly dependent random variables.
The Rosenblatt process 
can also be defined  as a double stochastic integral \cite{T4,Tu}.
Our approach  
consists in deriving the Rosenblatt
 process from a specific random 
particle 
system, which hopefully provides an intuitive physical interpretation of 
this
 process. A useful tool here is the theory of random variables in the 
space of tempered distributions ${\cal S}'\equiv{\cal S}'({\mathbb R}^d)$ 
($d=1$ in 
our case), since it permits to employ some nice properties of this space, 
and the Rosenblatt process can be expressed with the help of an ${\cal 
S}'$-random variable. This was noted by Dobrushin \cite{D}. Relations 
between random particle systems and random elements of ${\cal S}'$ have 
been 
studied by many authors, beginning with Martin-L\"of \cite{ML}. Our 
approach is in the spirit of Adler and coauthors, e.g. \cite{AE}, where a 
general scheme (but still not so general to cover our case) was developed 
for 
representing a given random element of ${\cal S}'$ as the 
limit of appropriate functionals of some particle system (with high 
density in \cite{AE}). That approach was later applied in \cite{AFL,BT1,Ta2}
to give  particle picture interpretations of the self-intersection local 
times of density processes in ${\cal S}'$. We stress that our principal 
aim is to construct the Rosenblatt process by means of a  
particle system, and to this end an important step is to study a suitable 
random element of ${\cal S}'$. Particle picture approaches have been used 
to obtain fractional Brownian motion and subfractional Brownian motion 
with Hurst parameter $H$, in different ways for $H<1/2$ and $H>1/2$ 
\cite{BGT1,BT2}.

Our results can be summarized as follows.

The Rosenblatt process with parameter $H$, defined for $H\in(1/2,1)$, is 
represented in the form 
$\xi_t=\langle Y,\UNO_{[0,t]}\rangle, t\geq 0$, where $Y$ is an 
${\cal S}'$-random variable which is obtained from the Wick product 
$:X\otimes 
X:$, where $X$ is a centered Gaussian ${\cal S}'$-random variable (see (\ref{eq:3.1})). $X$ is 
in a sense a distributional derivative of a suitable fractional Brownian motion, and 
the relation between $Y$ and $X$ corresponds to the informal formula 
(43) in \cite{T4}. Note that $:X\otimes X:$ is a counterpart of 
the Hermite polynomial of order $2, H_2(x)=x^2-1$, used in \cite{T1,T4}.  The 
possibility to take
 ``test 
functions'' of the form $\varphi =\UNO_{[0,t]}$ to bring in a time 
parameter has been noted for example in \cite{T2}. This formulation  is in 
the spirit 
of white noise analysis.

We define a particle system on ${\mathbb R}$ with initial distribution 
given by a Poisson random field with Lebesgue intensity measure, and 
particles evolving independently according to the standard symmetric 
$\alpha$-stable L\'evy process. The particles are independently assigned 
charges $+1$ and $-1$ with probabilities $1/2$. A crucial element of the 
construction is the intersection local time of two independent 
$\alpha$-stable processes, which is known to exist for $\alpha >1/2$ 
(Proposition 5.1 of \cite{BT1}). Consider the process $\xi^T$ defined by
\begin{equation}
\label{eq:1.1}
\xi^T_t=\frac{1}{T}\sum_{j,k\atop_{j\neq k}}\sigma_j\sigma_k\langle 
\Lambda(x_j+\rho^j,x_k+\rho^k;T), \UNO_{[0,t]}\rangle, \,\, t\geq 0,
\end{equation}
where the $x_j$ are the points of the initial Poisson configuration, the 
$\rho^j$ are the $\alpha$-stable processes corresponding to those points 
$(\rho^j(0)=0)$, the 
$\sigma_j$ are the respective charges, and 
$\Lambda(x_j+\rho^j,x_k+\rho^k;T)$ is the intersection local time of the 
processes $x_j+\rho^j$ and $x_k+\rho^k$ on the interval $[0,T]$. This 
local time is defined as a process in ${\cal S}'$ (Definition \ref{d:2.2} and Proposition \ref{p:2.3}). $\xi^T$ is 
given a precise meaning and shown to have a continuous modification (Lemma \ref{l:3.4}).
The result  is that for $\alpha\in (1/2,1),\xi^T$ converges 
as 
$T\to\infty$, in 
the sense of weak functional convergence, to the Rosenblatt process with 
parameter $H=\alpha$, up to a multiplicative constant (Theorem \ref{t:3.5}).

Our second result concerns the analysis of long range 
dependence of the 
Rosenblatt process.
 We give a precise measure of  long range dependence 
by means of a number called {\it dependence exponent}
 (Theorem \ref{t:4.1}). 
We show that the increments are asymptotically 
independent (not just uncorrelated) as the distance between the intervals 
tends to infinity (Corollary \ref{c:4.2}). 

In Section 2 we give background on the Rosenblatt process, on the ${\cal 
S}'$-random variable $X$ from which  fractional Brownian motion is 
derived, on how the abovementioned particle system produces $X$, and on 
intersection local time. In Section 3 we construct the Rosenblatt process 
by means of the ${\cal S}'$-random variable $Y$ obtained from $:X\otimes 
X:$, we show that the process $\xi^T$ is well defined, and we prove  
convergence to the Rosenblatt process. In Section 4 we discuss the long 
range dependence of the process. Section 5 contains some  
comments of related interest. The proofs are given in Section 6.

We use the following notation:\\
$\cal S$: the Schwartz space of smooth rapidly decreasing functions on 
$\mathbb R$,\\
 ${\cal S}'$: the space of tempered distributions (dual of ${\cal S}$), \\
$\overline{\phantom{\varphi}}$: complex conjugate,\\
$\widehat\varphi(z)=\int_{\mathbb R} e^{ixz}\varphi(x)dx$: Fourier transform of 
a function $\varphi$,\\
$\Rightarrow$: convergence in law in an appropriate space,\\
$\Rightarrow_f$: convergence of finite-dimensional distributions,\\
$*$: convolution,\\
$C([0,\tau])$: the space of real continuous functions on $[0,\tau]$,\\
$C, C_i$: generic positive constants with possible dependencies in parentheses. 

\section{Background}

\subsection{The Rosenblatt process}
\setcounter{equation}{0}

We recall some facts on the Rosenblatt process $\xi=(\xi_t)_{t\geq 0}$ with parameter $H\in (1/2,1)$, which are found in \cite{T1,T4}. The 
characteristic function of the finite-dimensional distributions of the process in a small neighborhood of 0 has the form
\begin{equation}
\label{eq:2.1}
E{\rm exp}\left\{i\sum ^p_{j=1}\theta_j\xi_{t_j}\right\}
={\rm 
exp}\left\{\frac{1}{2}\sum^\infty_{k=2}\frac{(2i\sigma)^k}{k}R_{H,k}(\theta_1,\ldots,\theta_p; t_1,\ldots,t_p)\right\},
\end{equation}
where
\begin{eqnarray}
\lefteqn{
R_{H,k}(\theta_1,\ldots,\theta_p;t_1,\ldots,t_p)
=\int_{{\mathbb R}^k}\psi(x_1)\psi(x_2)\ldots \psi(x_k)}\nonumber\\
\label{eq:2.2}
&&\cdot|x_1-x_2|^{H-1}|x_2-x_3|^{H-1}\ldots|x_{k-1}-x_k|^{H-1}|x_k-x_1|^{H-1}dx_1\ldots dx_k,
\end{eqnarray}

\begin{equation}
\label{eq:2.3}
\psi(x)=\sum^p_{j=1}\theta_j\UNO_{[0,t_j]}(x),
\end{equation}
and
\begin{equation}
\label{eq:2.4}
\sigma=\left[\frac{1}{2}H(2H-1)\right]^{1/2}.
\end{equation}
The value of $\sigma$ is chosen so that $E\xi^2_1=1$. The series in the exponent converges for $\theta_1,\ldots,\theta_p$ in a (small) neighborhood of $0$ depending on $t_1,\ldots, t_p$, and 
(\ref{eq:2.1}) defined in this neighborhood determines the distribution.
The process is also characterized by the cumulants of the random variable $\sum^p_{j=1}\theta_j\xi_{t_j}$, which are $\kappa_1=0$,
\begin{equation}
\label{eq:2.5}
\kappa_k=2^{k-1}(k-1)!\sigma^kR_{H,k}(\theta_1,\ldots,\theta_p;t_1,\ldots,t_p),\,\, 
k\geq 2.
\end{equation}

The process $\xi$ arises from a (Donsker-type) limit in distribution as $n\to\infty$ of the processes
$$\xi_n(t)=\frac{\sigma}{n^H}\sum^{\lfloor nt\rfloor}_{j=1}X_j,\,\,t\in [0,T],$$
($\lfloor \cdot\rfloor$ denotes integer part), where the random variables $X_j$ are defined by $X_j=Y^2_j-1$ (i.e., $X_j=H_2(Y_j)$, where $H_2$ is the Hermite polynomial of order 2), and  $(Y_j)_j$ is
a Gaussian stationary sequence of
  random variables with mean $0$, variance $1$, and covariances
$$r_j=EY_0Y_j=(1+j^2)^{(H-1)/2}\sim j^{H-1}\,\,{\rm as}\,\, j\to\infty.$$

The spectral  representation of the process is
\begin{equation}
\label{eq:2.6}
\xi_t\stackrel{d}{=}A(H)\int^{''}_{{\mathbb R}^2}\frac{e^{i(\lambda_1+\lambda_2)t}-1}{i(\lambda_1+\lambda_2)}\frac{1}{|\lambda_1|^{H/2}|\lambda_2|^{H/2}}d\widetilde{B}(\lambda_1)d\widetilde{B}(\lambda_2)
\end{equation}
($\stackrel{d}{=}$ means equality in distribution), where
\begin{equation}
\label{eq:2.7}
A(H)=\frac{[H(2H-1)/2]^{1/2}}{2\Gamma(1-H)\sin(H\pi/2)},
\end{equation}
$\widetilde{B}$ is a complex Gaussian measure on $\mathbb R$ such that 
$\widetilde{B}=B^{(1)}+iB^{(2)},B^{(1)}(A)=B^{(1)}(-A), 
B^{(2)}(A)\\=-B^{(2)}(-A), A$ is a Borel set of $\mathbb R$ with finite 
Lebesgue 
measure 
$|A|, B^{(1)}$ and $B^{(2)}$ are independent, and 
$E(B^{(1)}(A))^2=E(B^{(2)}(A))^2=\frac{1}{2}|A|$. ($\widetilde{B}$ can be viewed as a complex-valued Fourier transform of white noise). The double prime on the integral means that the diagonals  
$\lambda_1=\pm\lambda_2$ are excluded in the integration. The process also has a time  representation as a double integral on ${\mathbb R}^2$ with respect to Brownian motion, and a finite interval integral representation obtained in 
\cite{Tu}.

We have  mentioned in the Introduction some of the main properties of the 
Rosenblatt process. We recall the selfsimilarity with parameter $H$: 
for any $c>0$, 
\begin{equation}
\label{eq:2.8}
(\xi_{ct})_{t\geq 0}\stackrel{d}{=} c^H(\xi_t)_{t\geq 0}
\end{equation}
This
and stationarity of increments imply
\begin{equation}
\label{eq:2.9}
E(\xi_t-\xi_s)^2=\sigma^2(t-s)^{2H},
\end{equation}
hence the covariance function of $\xi$ has the same form as that of the fractional Brownian motion, i.e.,
\begin{equation}
\label{eq:2.10}
E\xi_s\xi_t=\frac{\sigma^2}{2}(s^{2H}+t^{2H}-|t-s|^{2H}).
\end{equation}
In particular,
 the increments are positively correlated (since $H>1/2$), 
and 
\begin{equation}
\label{eq:2.11}
\sum^\infty_{k=1}E\xi_1(\xi_k-\xi_{k-1})=\infty.
\end{equation}

We have not found a published proof of the non semimartingale property of 
$\xi$, but that is easy to show. By (\ref{eq:2.9})
 with $H>1/2$ it is obvious that  the quadratic variation is $0$. A deeper 
result is that  selfsimilarity and 
increment stationarity imply that the paths have infinite variation \cite{V}.
The non semimartingale property of fractional Brownian motion $(H\neq 
1/2)$ follows for example from a general criterion for Gaussian processes 
\cite{BGT0} (Lemma 2.1, Corollary 2.1).

Infinite divisibility was recently proved in \cite{VT,MT1}.

There does not seem
 to be information in the literature on whether or not the Rosenblatt 
process has the Markov property (it seems plausible that it does not, by 
analogy with fractional Brownian motion). 

The fractional Brownian motion is the only Gaussian process that has the properties (\ref{eq:2.8}) and (\ref{eq:2.9}) (with $\sigma=1$). On the
other hand, there are many processes with stationary increments satisfying 
(\ref{eq:2.8}) which belong to the second chaos \cite{MT}. The Rosenblatt 
process is the simplest one of them.

\subsection{An ${\cal S}'$-random variable related to fractional Brownian motion}

Recall that fractional Brownian motion (fBm) with Hurst parameter 
$H\in(0,1)$ is a centered Gaussian process $(B^H_t)_{t\geq 0}$ with 
covariance given by the right hand side of 
(\ref{eq:2.10}) with $\sigma=1$ (see \cite{ST} for background on fBm).  
This process can be represented with the help of the centered Gaussian ${\cal S}'$-valued random variable $X$ with covariance functional
\begin{equation}
\label{eq:2.12}
E\langle X,\varphi_1\rangle \langle X,\varphi_2\rangle=\frac{1}{\pi}
\int_{\mathbb R}\widehat{\varphi}_1(x)\overline{\widehat{\varphi}_2(x)}|x|^{-\alpha}dx,\quad \varphi_1,\varphi_2\in{\cal S},
\end{equation}
where $-1<\alpha<1$. Namely, 
\begin{equation}
\label{eq:2.13}
\left(\langle X,\UNO_{[0,t]}\rangle\right)_{t\geq 0}=(KB^H_t)_{t\geq 0},
\end{equation}
with $H=\frac{1+\alpha}{2}$ and $K$ is some positive constant. ($\langle X,\UNO_{[0,t]}\rangle$ is defined by an $L^2$-extension). For $0<\alpha<1\quad (\frac{1}{2}< H<1)$ this random variable $X$ can be obtained from the particle system described in the Introduction, i.e., the system of independent standard $\alpha$-stable processes (particle motions) starting from a Poisson random field with Lebesgue
 intensity. Each particle has a charge $\pm1$ with equal probabilities, and the charges are mutually independent and independent of the initial configuration and of the particle motions. The motions have the form $x_j+\rho^j$, where the $x_j$'s are the points of the intial configuration, the $\rho^j$ are independent standard $\alpha$-stable processes independent of $\{x_j\}_j, \rho^j_0=0$. The charges are denoted by $\sigma_j$.

The normalized total charge occupation on the  interval $[0,T]$ is defined by
\begin{equation}
\label{eq:2.14}
\langle X_T,\varphi\rangle=\frac{1}{\sqrt{T}}\sum_j\sigma_j\int^T_0\varphi(x_j+\rho^j_s)ds,\quad\varphi\in {\cal S}.
\end{equation}

We have the following proposition.

\begin{prop}
\label{p:2.1}
If $0<\alpha<1$ and $T\to\infty$, then 

\noindent
(a) $X_T\Rightarrow X$ in ${\cal S}'$, where $X$ is as in (\ref{eq:2.12}).

\noindent
(b) $(\langle X_T,\UNO_{[0,t]}\rangle)_{t\geq 0}\;{\Rightarrow}_f
(KB^H_t)_{t\geq 0}$ with $H=\frac{1+\alpha}{2}$.
\end{prop}

This fact is an easy consequence of Theorem 2.1(a) in \cite{BGT2}, where the occupation time fluctuations around the mean for the system without charges were considered. It suffices to take two independent copies of such systems and to write the difference of their occupation time fluctuations.

A similar procedure with a different functional of a particle system without charges permits also to obtain fBm with $H<\frac{1}{2}$, as well as the corresponding random variable $X$ (see Theorem 2.9 in \cite{BT2}).

\subsection{Intersection local time}
\label{sub:2.3}

There are several ways to define the intersection local time (ILT) of two processes (see e.g. \cite{Dy,Ro,AE}. We will take the definition from \cite{BT1}, which is close to that of \cite{AE}. Intuitively, ILT of real cadlag processes 
$(\eta^1_t)_{t\geq 0},(\eta^2_t)_{t\geq 0}$ up to time $T$ is given by
$$\langle\Lambda(\eta^1,\eta^2;T),\varphi\rangle=\int^T_0\int^T_0
\delta(\eta^2_v-\eta^1_u)\varphi(\eta^1_u)dudv,\quad\varphi\in{\cal S},$$
where $\delta$ is the Dirac distribution. We want to regard $\Lambda$ as a process in ${\cal S}'$.

To make this definition rigorous one has to apply a limiting procedure.

Let ${\cal F}$ denote the class of nonnegative symmetric infinitely differentiable functions $f$ on $\mathbb R$ with compact support and such that 
$\int_{\mathbb R} f(x)dx=1$. For $f\in {\cal F}, \varepsilon >0$, let 
\begin{equation}
\label{eq:2.15}
f_\varepsilon(x)=\varepsilon^{-1}f \left(\frac{x}{\varepsilon}\right),
\quad x\in{\mathbb R}.
\end{equation}
We will frequently use 
\begin{equation}
\label{eq:2.16}
|\widehat{f}_\varepsilon{(x)}|\leq 1\quad{\rm and}\quad 
\lim_{\varepsilon\to 0}\widehat{f}_\varepsilon(x)=\lim_{\varepsilon\to 0}\widehat{f}(\varepsilon x)=1.
\end{equation}
We define
\begin{equation}
\label{eq:2.17}
\langle\Lambda^f_\varepsilon(\eta^1,\eta^2;T),\varphi\rangle=
\int^T_0
\int^T_0f_\varepsilon(\eta^2_v-\eta^1_u)\varphi(\eta^1_u) dudv,\quad T\geq 0,\varphi\in{\cal S}.
\end{equation}

\begin{defn}
\label{d:2.2}
If there exists an ${\cal S}'$-process
$\Lambda(\eta^1,\eta^2)=(\Lambda(\eta^1,\eta^2;T))_{T\geq 0}$ such that for each $T\geq 0, \varphi\in{\cal S}$ and any $f\in{\cal F},\langle\Lambda(\eta^1,\eta^2;T),\varphi\rangle$ is the mean square limit of $\langle\Lambda^f_\varepsilon(\eta^1,\eta^2;T),\varphi\rangle$ as $\varepsilon\to 0$, then the process $\Lambda(\eta^1,\eta^2)$ is called the intersection local time (ILT) of the processes $\eta^1,\eta^2$.
\end{defn}

In \cite{BT1} the following result was proved (Theorem 4.2 and Proposition 5.1 therein).

\begin{prop}
\label{p:2.3}

Let $\eta^1,\eta^2$ be independent standard $\alpha$-stable processes in $\mathbb R$. If $\alpha>\frac{1}{2}$, then for any $x,y\in\mathbb R$ the ILT $\Lambda(x+\eta^1,y+\eta^2)$ exists. Moreover, for all $T\geq 0, f\in {\cal F},
 \varphi\in {\cal S},\\ \langle\Lambda^f_\varepsilon(\cdot +\eta^1,\cdot +\eta^2;T),\varphi\rangle$ converges in $L^2(\mathbb R^2\times \Omega,\lambda\otimes\lambda\otimes P)$, where $\lambda$ is the Lebesgue measure on $\mathbb R$, and $P$ is the probability measure on the underlying sample space $\Omega$.
\end{prop}

\section{Particle picture for the Rosenblatt process}
\setcounter{section}{3}
\setcounter{equation}{0}

We begin with another representation of the Rosenblatt process, which is more suitable for our purpose. From \cite{D} it can be deduced that this construction was known to Dobrushin, but we have not been able to find it written explicitly in the literature. Therefore we will describe it in detail, but the proof will be only sketched.

Let $X$ be the centered Gaussian ${\cal S}'$-random variable with covariance (\ref{eq:2.12}). Recall that the Wick product $:X\otimes X:$ is defined as a random variable in ${\cal S}'(\mathbb R^2)$ such that
\begin{equation}
\label{eq:3.1}
\langle :X\otimes X:,\varphi_1\otimes\varphi_2\rangle=\langle X,\varphi_1\rangle\langle X,\varphi_2\rangle-E\langle X,\varphi_1\rangle\langle X,
\varphi_2\rangle,\quad \varphi_1,\varphi_2\in{\cal S}.
\end{equation}
(see e.g. Chapter 6 of \cite{GJ} or \cite{AR,AT}). 
 The Wick square of $X$ is an ${\cal S}'$-random variable $Y$ that  can be written informally as $\langle Y,\varphi\rangle=\langle :X\otimes X:, \varphi(x)\delta_{y-x}\rangle$. To make this rigorous we use approximation. Fix $f\in{\cal F}$ and let $f_\varepsilon$ be as in (\ref{eq:2.15}). For $\varphi\in{\cal S}$ we denote
\begin{equation}
\label{eq:3.2}
\Phi^f_{\varepsilon,\varphi}(x,y)=\varphi(x)f_\varepsilon(y-x),
\end{equation}
and we define an ${\cal S}'$-random variable $Y^f_\varepsilon$ by
\begin{equation}
\label{eq:3.3}
\langle Y^f_\varepsilon,\varphi\rangle=\langle :X\otimes X:,\Phi^f_{\varepsilon,\varphi}\rangle,\quad \varphi\in{\cal S}.
\end{equation}
The following lemma is an easy consequence of the regularization theorem 
\cite{I}, and the fact that, by Gaussianity,
\begin{equation}
\label{eq:3.4}
E\langle :X\otimes X:,\Phi\rangle\langle :X\otimes X:,\Psi\rangle=\frac{1}{\pi^2} \int_{\mathbb R}\widehat{\Phi}(x,y)(\overline{\widehat{\Psi}(x,y)}+\overline{\widehat{\Psi}(y,x)})|x|^{-\alpha}|y|^{-\alpha}dxdy,\quad\Phi,\Psi\in{\cal S}(\mathbb R^2).
\end{equation}

\begin{lem}
\label{l:3.1}
If $\frac{1}{2}<\alpha<1$, then there exists an ${\cal S}'$-random variable $Y$ such that for any $f\in{\cal F}$,
$$\langle Y,\varphi\rangle=L^2 {\rm -} \lim_{\varepsilon\to 0}\langle Y^f_\varepsilon,\varphi\rangle,\quad \varphi\in{\cal S}.$$%
$\langle Y,\varphi\rangle$ can be further extended in $L^2(\Omega)$ to test functions of the form $\UNO_{[0,t]}$.
\end{lem}

The next theorem is an analogue of (\ref{eq:2.13}) for the Rosenblatt process.

\begin{thm}
\label{p:3.2}
Let $Y$ be as in Lemma \ref{l:3.1}. Then the real process
\begin{equation}
\label{eq:3.5}
(\langle Y,\UNO_{[0,t]}\rangle)_{t\geq 0}
\end{equation}
is, up to a constant, the Rosenblatt process with parameter $H=\alpha$.
\end{thm}
\vglue .25cm

We remark that this proposition gives a rigorous meaning to the informal expression (43) in \cite{T4} relating the Rosenblatt process $\xi$ with parameter $H$ and a fBm with parameter
 $H_1 =\frac{H+1}{2}\in (\frac{3}{4},1)$, which is given by
$$
\xi_t =C(H_1)\int^t_0 \left((B^{H_1}_s)'\right)^{2}ds.
$$
$X$ corresponds to $(B^{H_1})'$ and $:X\otimes X:$ 
corresponds to $((B^{H_1})')^2$. 
The relationship between the parameters follows from Proposition \ref{p:2.1}(b) and Theorem \ref{p:3.2}.

In \cite{D}, ${\cal S}'$-random variables such as $Y$ are represented in terms of complex multiple stochastic integrals related to (\ref{eq:2.6}).

After the first version of this paper had been submitted, the referee drew our attention to  preprint \cite{Arras} which appeared in the meantime. That paper uses Hida-Kuo  type calculus \cite{Kuo} to construct the stochastic integral with respect to the Rosenblatt process, but the representation \eqref{eq:3.5} does not seem to be present there.

Representation (\ref{eq:3.5}) and Proposition \ref{p:2.1} suggest a way to construct the Rosenblatt process by means of a particle system. We consider the particle system as before with $\frac{1}{2}<\alpha<1$. By Proposition \ref{p:2.3}, for each pair $\rho^j,\rho^k, j\neq k$, the intersection local time $\Lambda(x_j+\rho^j,x_k+\rho^k;T)$ exists, moreover it extends in a natural way to test function $\UNO_{[0,t]}$. Namely, we have the following lemma.

\begin{lem}
\label{l:3.3}
Let
\begin{equation}
\label{eq:3.6}
\psi=\sum^m_{j=1}a_j\UNO_{I_j},\quad a_j\in\mathbb R,\,\, 
I_j\,\,{\it is\,\, a}\,\, {\it bounded\,\, interval},
\end{equation}
Fix $\rho^1,\rho^2$, independent standard $\alpha$-stable processes and $x,y\in\mathbb R, T>0,f\in {\cal F}$.

\begin{itemize}
\item [(a)]
There exists an $L^2$-limit of
$\langle\Lambda^f_\varepsilon(x+\rho^1,y+\rho^2;T),\psi\rangle\quad{\it as}\quad \varepsilon\to 0$,
where $\Lambda^f_\varepsilon$ is given by (\ref{eq:2.17}), and this limit does not depend on $f$. We denote it by $\langle\Lambda(x+\rho^1,y+\rho^2;T),\psi\rangle$.
\item[(b)]
$$\langle\Lambda(x+\rho^1,y+\rho^2;T),\psi\rangle=L^2-\lim_{\varepsilon\to 0}\langle \Lambda(x+\rho^1,y+\rho^2;T),\psi *f_\varepsilon\rangle.$$
\item[(c)]
Moreover, $\langle\Lambda^f_\varepsilon(\cdot +\rho^1,\cdot +\rho^2;T),\psi\rangle$ converges in $L^2(\mathbb R^2\times\Omega,\lambda\otimes\lambda\otimes P)$.
\end{itemize}
\end{lem}

For the convenience of the reader let us recall \eqref{eq:1.1},
\begin{equation*}
\xi^T_t=\frac{1}{T}\sum_{j,k\atop_{j\neq k}}\sigma_j\sigma_k\langle 
\Lambda(x_j+\rho^j,x_k+\rho^k;T), \UNO_{[0,t]}\rangle, \,\, t\geq 0.
\end{equation*}

\begin{lem}
\label{l:3.4}
The process $\xi^T$ 
is well defined (the series converges in $L^2$), and it has a continuous modification.
\end{lem}

The main result of the paper is stated in the next theorem, which is a counterpart of Proposition \ref{p:2.1}.

\begin{thm}
\label{t:3.5}
Let $\frac{1}{2}<\alpha<1$. Then $\xi^T\Rightarrow K\xi$ in $C([0,\tau])$ as 
$T\to\infty, \tau>0$, where $\xi$ is the Rosenblatt process with $H=\alpha$ and $K$ is a positive constant.
\end{thm}

\section{Dependence exponent}
\setcounter{section}{4}
\setcounter{equation}{0}

Long range dependence is a general notion that has not been clearly defined and can be viewed in different ways \cite{S,T3,H}. For a Gaussian process $\eta$, long range dependence is usually described as slow (power) decay of the covariance of increments on intervals $[u,v], [s+\tau, t+\tau]$ as $\tau\to\infty$, i.e.,
$${\rm Cov}(\eta_v-\eta_u,\eta_{t+\tau}-\eta_{s+\tau})\sim C_{u,v,s,t}\tau^{-K},$$
where $K$ is a positive constant, and convergence or divergence of the series
$$\sum^\infty_{k=1} {\rm Cov}(\eta_1-\eta_0, \eta_{k+1}-\eta_k)$$
are sometimes referred to as ``short range'' dependence and ``long range'' dependence, respectively. This criterion is also applied to non Gaussian processes with finite second moments, such as the Rosenblatt process \cite{T4}. The underlying idea is that the increments become uncorrelated (but not necessarily independent) at some rate as the distance $\tau$ between the intervals tends to infinity. However, it can happen that $K\leq 0$, which should also be regarded as long range dependence (\cite{G} contains examples).

In order to characterize long range dependence in some more precise way for infinitely divisible processes (not necessarily Gaussian), the codifference (see \cite{RZ,ST}) can be useful. In \cite{BGT3} we defined the {\it dependence exponent} of a (real) infinitely divisible process $\eta$ as the number 
\begin{equation}
\label{eq:4.1}
\kappa=\inf_{z_1,z_2\in{\mathbb R}}\inf_{0\leq u<v<s<t}
\sup\{\gamma>0:D^\eta_\tau(z_1,z_2;u,v,s,t)=o(\tau^{-\gamma})\,\,{\rm as}\,\, 
\tau\to\infty\},
\end{equation}
where
\begin{eqnarray}
D^\eta_\tau(z_1,z_2;u,v,s,t)&=&|\log Ee^{i(z_1(\eta_v-\eta_u)+z_2(\eta_{t+\tau}-\eta_{s+\tau})}\nonumber\\
\label{eq:4.2}
&&-\log Ee^{iz_1(\eta_v-\eta_u)}-\log Ee^{iz_2(\eta_{t+\tau}-\eta_{s+\tau})}|
\end{eqnarray}
is the absolute value of the codifference of the random variables $z_1(\eta_v-\eta_u)$ and $-z_2(\eta_{t+\tau}-\eta_{s+\tau})$. Note that if $\eta$ is Gaussian, then
$$D^\eta_\tau(z_1,z_2;u,v,s,t)=|z_1z_2{\rm Cov}(\eta_v-\eta_u,\eta_{t+\tau}-\eta_{s+\tau})|.$$
For fractional Brownian motion, $\kappa=K=2-2H$, and for sub-fractional Brownian motion, $\kappa=K=3-2H$ \cite{BGT1}.

It turns out that the same idea can be used to measure long range 
dependence for the 
Rosenblatt process, and this can be done without recourse to infinite divisibility. 
As recalled in Subsection 2.1, the characteristic functions of the 
finite-dimensional distributions of the process are given by an explicit 
formula only for small values of the parameters, which are $z_1$ and $z_2$ 
in our case (see (\ref{eq:4.2})). We show next that it is enough to take 
$z_1$ and 
$z_2$ in an appropriate neighborhood of $0$ to measure long range 
dependence and  prove 
asymptotic independence of increments. 

For simplicity we take $u=s,v=t$.

\begin{thm}
\label{t:4.1}
Let $\xi$ be the Rosenblatt process with parameter $H$. For any $0\leq s<t$ there exists a neighborhood $U(s,t)$ of $0$ in $\mathbb R^2$ such that
\begin{equation}
\label{eq:4.3}
 D^\xi_\tau(z_1,z_2,s,t):= {D}^\xi_\tau(z_1,z_2,s,t,s,t)
\end{equation}
is well defined for $(z_1,z_2)\in U(s,t)$ and all $\tau>0$, and if we modify 
(\ref{eq:4.1}) putting
\begin{equation}
\label{eq:4.4}
\kappa=\inf_{0\leq s<t}\inf_{(z_1,z_2)\in U(s,t)}\sup\{\gamma>0:{D}^\xi_\tau(z_1,z_2,s,t)=o(\tau^{-\gamma})\quad{\rm as}\quad \tau\to\infty\},
\end{equation}
then $\kappa=2-2H$.
\end{thm}

So we see that dependence exponent of the Rosenblatt process with parameter $H$ is the same as that for fBm $B^H$.

From this theorem, by a standard tightness argument, stationarity of increments of $\xi$, and the fact that the law of $\xi_t$ is determined by its characteristic function in an arbitrarily small neighborhood of $0$, we obtain the following corollary.

\begin{cor}
\label{c:4.2}
For any $0<s<t$, the increments of the Rosenblatt process $\xi_t-\xi_s$ and $\xi_{t+\tau}-\xi_{s+\tau}$ are asymptotically independent as $\tau\to\infty$, i.e., if $\mu_{(s,t)}$ is the law of $\xi_t-\xi_s$ and $\mu_{(s,t),(s+\tau,t+\tau)}$ is the law of $(\xi_t-\xi_s,\xi_{t+\tau}-\xi_{s+\tau})$, then
$$\mu_{(s,t),(s+\tau,t+\tau)}\Rightarrow\mu_{(s,t)}\otimes\mu_{(s,t)}\left(=\mu_{(s,t)}\otimes\mu_{(s+\tau,t+\tau)}\right).$$
\end{cor}

\section{Additional comments}
\setcounter{section}{5}
\setcounter{equation}{0}

\subsection{Sub-Rosenblatt process}
\label{sub:5.1}

It is known that if in the formula (\ref{eq:2.13}) we put $\UNO_{[0,t]}-\UNO_{[-t,0]}$ instead of $\UNO_{[0,t]}$, we obtain a sub-fractional Brownian motion 
(sub-fBm), i.e., a centered Gaussian process with covariance 
$$t^{2H}+s^{2H}-\frac{1}{2}(|t-s|^{2H}+(t+s)^{2H}),$$
again with $H=\frac{1+\alpha}{2}$. This process has been studied by 
several 
authors, e.g., \cite{BGT1,DZ,T,YSH}
and others. In particular, in \cite{BT2} an analogue of Proposition \ref{p:2.1}(b) was proved for sub-fBm.

We can now extend  formula (\ref{eq:3.5}), and define a new process $(\langle Y,\UNO_{[0,t]}-\UNO_{[-t,0]}\rangle)_{t\geq 0}$. It is natural to call it sub-Rosenblatt process, as it has the same covariance as sub-fBm.

Analogues of  Theorems \ref{t:3.5} and \ref{t:4.1} also hold.

\subsection{Rosenblatt process with two parameters}
\label{sub:5.2}

Maejima and Tudor in \cite{MT} define a class of self-similar processes with stationary increments that live in the second Wiener chaos. These processes depend on two parameters $H_1,H_2$, and the Rosenblatt process corresponds to the case $H_1=H_2$. One can ask about the possibility of extending our construction to those two-parameter processes.

\subsection{General Hermite processes}
\label{sub:5.3}

Taqqu \cite{T2} studies extensions of the Rosenblatt process living in Wiener chaos of order $k, k\geq 2$, which he calls Hermite processes. One can attempt to find a particle picture interpretation for those processes. It seems that one should employ $k$-th Wick powers and work with intersection local times of 
$k$-tuples of stable processes.

\section{Proofs}
\setcounter{section}{6}
\setcounter{equation}{0}

\noindent 
{\bf Proof of Theorem \ref{p:3.2} (outline)}
Let $\xi$ be the Rosenblatt process with parameter $H$. It is known that its distributions are determined by its moments, therefore it is enough to prove that for all $n,p\in \ene$, $t_1,t_2, \ldots, t_p\ge 0$,
\begin{equation}
 \label{e:p.1}
 E\left <Y,\psi\right >^n=C^nE \left(\sum_{j=1}^p\theta_j\xi_{t_j}\right)^n,
\end{equation}
where $\psi$ has the form \eqref{eq:2.3}. It is known (see e.g. \cite{Shiryaev}, Thm II.12.6) that
\begin{equation}
 \label{e:p.2}
 E \left(\sum_{j=1}^p\theta_j\xi_{t_j}\right)^n=\sum_{\pi\in \Pn}\prod_{B\in \pi}\kappa_{\# B},
\end{equation}
where $\kappa$'s are the corresponding cumulants given by \eqref{eq:2.5} and $\Pn$ is the set of all partitions of $\{1,\ldots, n\}$, and $\#$ denotes cardinality of a set. 

To compute $E\left <Y,\psi\right >^n$  we will need the formulas for moments of the Wick product $:X\otimes X:$. These moments are expressed with the help of Feynmann graphs (see e.g. \cite{AR} p.\ 1085 or \cite{Ta1}, p. 422). For fixed $n$ we consider graphs as follows. Suppose that 
 we have  
$n$ numbered vertices.
Each vertex has two legs numbered $1$ and $2$.
Legs are paired, forming links between vertices, in such a way that
each link connects two different vertices, and there are no unpaired legs left.
The graph is a set of links. Each link is described by an unordered pair
$\{(i,j),(l,m)\}$, $i,l\in\{1,\ldots n\}$, $j,m\in\{1,2\}$,
which means that leg  $j$, growing from vertex $i$ is paired with leg
$m$ growing from vertex $l$. $i\neq l$ since each link connects different
vertices, and any $(i,j)$, $i=1,\ldots ,n$, $j=1,2$, is a part of one
and only one link. The set of all distinct  graphs of the above form will be
denoted by $\G_n^2$. 
Let $\tilde \G_n^2$ denote the set of all connected graphs in $\G_n^2$.

By formulas (2.7) and (6.1) in \cite{Ta1} we have 
\begin{equation}
 \label{e:p.3}
 E \left< :X\otimes X:, \Phi\right>^n=\sum_{G\in \G_n^2}I^G(\Phi), \qquad \Phi \in \Ss(\R^{2d}),
\end{equation}
where 
\begin{equation}
\label{e:p.4}
 I^G(\Phi)=\int_{\R^{2n}}\hat\Phi(x_{1,1},x_{1,2})\ldots \hat\Phi(x_{n,1}, x_{n,2})\prod_{\{(l,m),(p,q)\}\in G}\delta_{-x_{p,q}}(dx_{l,m})\abs{x_{p,q}}^{-H}dx_{p,q}.
\end{equation}
Using Lemma \ref{l:3.1} and similar arguments as in the proof of Lemma \ref{l:3.3} below it is not difficult to see that 
\begin{equation}
 \label{e:p.4a}
  E \left<Y, \psi\right>^n=\lim_{\varepsilon\to 0}\sum_{G\in \G_n^2}I^G(\Phi^f_{\varepsilon,\psi}),
\end{equation}
where $\Phi^f_{\varepsilon,\psi}$ is given by \eqref{eq:3.2}. 

For $G\in \tilde \G_n^2$ we have 
\begin{equation*}
 I^G(\Phi^f_{\varepsilon,\psi})=\int_{\R^{n}}\hat \psi(x_1-x_2)\hat\psi(x_2-x_3)\ldots \hat \psi(x_n-x_1) F^G_{\varepsilon}(x_1,\ldots, x_n)\abs{x_1}^{-H}\ldots \abs{x_n}^{-H}dx_1\ldots dx_n,
\end{equation*}
where $F^G_\varepsilon$ is a product of functions of the form $\hat f_\varepsilon(x_i)$ or $\overline{\hat f_\varepsilon(x_i)}$. By \eqref{eq:2.16} and the dominated convergence theorem,
\begin{equation}
 \label{e:p.5}
 \lim_{\varepsilon\to 0}I^G(\Phi^f_{\varepsilon,\psi})=J_n:=\int_{\R^n}\hat \psi(x_1-x_2)\hat\psi(x_2-x_3)\ldots \hat \psi(x_n-x_1) \abs{x_1}^{-H}\ldots \abs{x_n}^{-H}dx_1\ldots dx_n.
\end{equation}
In the proof of the integrability of the function under the integral in \eqref{e:p.5} we use \\
$|\hat \psi(x_1-x_2)\hat\psi(x_n-x_1)|\le |\hat \psi(x_1-x_2)|^2+|\hat \psi(x_n-x_1)|^2$,
\begin{equation}
\label{eq:6.3}
|\widehat{\psi}(x)|\leq \frac{C}{1+|x|},
\end{equation}
and 
\begin{equation*}
 \int_{\R}\frac{1}{1+\abs{x-y}}\abs{y}^{-H}dy\le C_1.
\end{equation*}

In the general case, if $G\in \G^2_n$, then it has a decomposition of the form $G=G_1\cup\ldots \cup G_k$, where $G_j$ are connected components of $G$. Then by \eqref{e:p.4} and \eqref{e:p.5} we have
\begin{equation}
 \label{e:p.7}
 \lim_{\varepsilon\to 0}I^G(\Phi^f_{\varepsilon,\psi})=J_{\# G_1}\ldots J_{\# G_k}, 
\end{equation}
where $\# G_j$ is the number of vertices of $G_j$.

Note that each $G\in \G_n^2$ determines a partition $\pi_G\in \Pn$, whose elements are the sets of vertices of connected components of $G$. Hence by \eqref{e:p.4a} and \eqref{e:p.7},
\begin{equation*}
 E\left<Y,\psi\right>^n=\sum_{G\in \G_n^2}\prod_{B\in \pi_G}J_{\# B}.
\end{equation*}
It is not difficult to see that if $\pi\in \Pn$ is of the form $\pi=\{B_1, \ldots, B_k\}$, $B_j\ge 2$, $j=1,\ldots, k$, then the number of different $G\in \G_n^2$ such that $\pi_G=\pi$ is equal to
\begin{equation*}
 \prod_{j=1}^k 2^{\# B_j-1}(\# B_j-1)!
\end{equation*}
Therefore, setting $J_1=0$ we obtain
\begin{equation}
 \label{e:p.8}
  E\left<Y,\psi\right>^n= \sum_{\pi\in \Pn}\prod_{B\in \pi} 2^{\# B-1}(\# B-1)! J_{\# B}.
\end{equation}
$J_k$ given by \eqref{e:p.5} can be also written as 
\begin{equation*}
 J_k=C^k\int_{\R^k}\psi(x_1)\ldots \psi(x_k)\abs{x_1-x_2}^{H-1}\abs{x_2-x_3}^{H-1}\ldots \abs{x_k-x_1}^{H-1}dx_1\ldots dx_k,
\end{equation*}
hence combining \eqref{e:p.8} with \eqref{eq:2.2}, \eqref{eq:2.5} and \eqref{e:p.2} we obtain \eqref{e:p.1}.
\qed

\bigskip

\noindent
{\bf Proof of Lemma \ref{l:3.3}} To prove part (a) it suffices to show that for any 
$f, g\in {\cal F}$,
\begin{equation}
\label{eq:6.1}
E\langle\Lambda^f_\varepsilon(x+\rho^1,y+\rho^2;T), \psi\rangle\langle\Lambda^g_\delta(x+\rho^1,y+\rho^2;T),\psi\rangle
\end{equation}
has a finite limit as $\varepsilon,\delta\to 0$. Analogously as in 
(7.1)-(7.4) of \cite{BT1},
 (\ref{eq:6.1}) is equal to
\begin{eqnarray*}
\lefteqn{
\frac{1}{(2\pi)^4}\int_{[0,T]^4}\int_{\mathbb R^4}e^{-ix(z+z')}e^{-iy(w+w')}\widehat{\psi}(z+w)\widehat{\psi}(z'+w')\widehat{f}(\varepsilon w)\widehat{g}(\delta w')}\\
&&\kern2.2cm\cdot\,\,\overline{\widehat{\nu}_{s,u}(z,z')}\overline{\widehat{\nu}_{r,v}(w,w')}dzdz'dwdw'dsdudrdv,
\end{eqnarray*}
where $\nu_{s,u}$ is the law of $(\rho^1_s,\rho^1_u)$. To complete the proof of part (a) it suffices to show that
\begin{equation}
\label{eq:6.2}
I:=\int_{\mathbb R^4}|\widehat{\psi}(z+w)||\widehat{\psi}(z'+w')|\int_{[0,T]^2}|\widehat{\nu}_{s,u}(z,z')|dsdu\int_{[0,T]^2}|\widehat{\nu}_{r,v}(w,w')|drdv
dzdz'dwdw'<\infty.
\end{equation}

To derive this we cannot repeat the argument of \cite{BT1} because $\widehat{\psi}\notin L^1$ for $\psi$ of the form (\ref{eq:3.6}), we only have \eqref{eq:6.3}.

Fix $\gamma>0$ such that $\frac{1}{2}+4\gamma<\alpha$. It is easy to see that
\begin{eqnarray}
\int_{[0,T]^2}|\widehat{\nu}_{s,u}(z,z')|dsdu&\leq& \frac{C(T)}{1+|z+z'|^\alpha}
\left(\frac{1}{1+|z|^\alpha}+\frac{1}{1+|z'|^\alpha}\right)\nonumber\\
\label{eq:6.4}
&\leq&C_1(T)h_\gamma(z,z')\frac{1}{1+|z|^\gamma}\frac{1}{1+|z'|^\gamma},
\end{eqnarray}
where
$$h_\gamma(z,z')=\frac{1}{1+|z+z'|^{\frac{1}{2}+\gamma}}\left(\frac{1}{1+|z|^{\frac{1}{2}+\gamma}}+\frac{1}{1+|z'|^{\frac{1}{2}+\gamma}}\right).$$
We have used
\begin{equation}
\label{eq:6.5}
\frac{1}{1+|z+z'|^\gamma}\leq C\frac{1+|z|^\gamma}{1+|z'|^\gamma}.
\end{equation}
Using $h_\gamma(z,z')h_\gamma(w,w')\leq h^2_\gamma(z,z')+h^2_\gamma(w,w'),$
(\ref{eq:6.4}) and symmetry, we obtain
$$I\leq C_2(T)\int_{\mathbb R^4}\frac{|\widehat{\psi}(z+w)|\,|\widehat{\psi}(z'+w')|}{(1+|z|^\gamma)(1+|w|^\gamma)\, (1+|z'|^\gamma)(1+|w'|^\gamma)}h^2_\gamma(z,z')dzdz'dwdw',$$
and (\ref{eq:6.5}) permits to replace the denominator by $(1+|z+w|^\gamma)(1+|z'+w'|^\gamma)$, hence (\ref{eq:6.2}) follows by (\ref{eq:6.3}).

Note that we have also shown that
\begin{eqnarray}
\lefteqn{
E\langle\Lambda(x+\rho^1,y+\rho^2;T),\psi\rangle^2
=\frac{1}{(2\pi)^4}\int_{[0,T]^4}\int_{\mathbb R^4}e^{-ix(z+z')}e^{-iy(w+w')}}\nonumber\\
\label{eq:6.6}
&&\kern1cm \cdot\widehat{\psi}(z+w)\widehat{\psi}(z'+w')
\overline{\widehat{\nu}_{s,u}(z,z')}\overline{\widehat{\nu}_{r,v}(w,w')}
dzdz'dwdw'dsdudrdv.
\end{eqnarray}

To prove part (b) we observe that the argument above can be carried out for linear combinations of functions of the form (\ref{eq:3.6}) and from ${\cal S}$ instead of $\psi$. Hence we see that (\ref{eq:6.6}) holds for $\psi-\psi * f_\varepsilon$, since 
\begin{equation}
\label{eq:6.7}
(\psi-\psi * f_\varepsilon)^{\widehat{}}(x)=\widehat{\psi}(x)(1-\widehat{f}(\varepsilon x))
\end{equation}

Then (b) follows from (\ref{eq:6.2}).

The proof of part (c) is the same as that of Proposition 4.4 in \cite{BT1}. Only the fact that $\psi\in L^2$ is needed here. $\hfill\Box$

\begin{rem}
\label{r:6.1}
{\rm From the proof of part (c) it follows that}
\begin{eqnarray}
\lefteqn{
E\int_{\mathbb R^2}(\langle\Lambda_\varepsilon^f(x+\rho^1,y+\rho^2;T),\psi\rangle-\langle\Lambda(x+\rho^1,y+\rho^2;T),\psi\rangle)^2dxdy}\nonumber\\
\label{eq:6.8}
&=&\frac{1}{(2\pi)^2}\int_{[0,T]^4}\int_{\mathbb R^2}|\widehat{\psi}(x+y)|^2|\widehat{f}(\varepsilon y)-1|^2e^{-|s-u||x|^\alpha}e^{-|r-v||y|^\alpha}dxdydsdrdudv,
\end{eqnarray}
{\rm and}
\begin{eqnarray}
\lefteqn{
E\int_{\mathbb R^2}(\langle\Lambda(x+\rho^1,y+\rho^2;T),\psi\rangle^2dxdy}
\nonumber\\
\label{eq:6.9}
&=&\frac{1}{(2\pi)^2}\int_{[0,T]^4}\int_{\mathbb R^2}|\widehat{\psi}(x+y)|^2|
e^{-|s-u||x|^\alpha}e^{-|r-v||y|^\alpha}dxdydsdrdudv.
\end{eqnarray}
\end{rem}

We need the following lemma which can be proved by repeating the argument of the proof of Lemma 4.1 in \cite{BT1}

\begin{lem}
\label{l:6.2}
For any $F(\cdot+\rho^1,\cdot +\rho^2)\in L^2(\mathbb R^2\times 
\Omega,\lambda\otimes\lambda\otimes P)$ the series
$$\sum_{j,k\atop j\neq k}\sigma_j\sigma_k F(x_j+\rho^j,x_k+\rho^k)$$
converges in $L^2(\Omega)$, and 
\begin{eqnarray}
\lefteqn{
E\left(\sum_{j,k\atop j\notin k}\sigma_j\sigma_k F(x_j+\rho^j,x_k+\rho^k)\right)^2}\nonumber\\
\label{eq:6.10}
&=&\int_{\mathbb R^2}E(F^2(x+\rho^1,y+\rho^2)+F(x+\rho^1,y+\rho^2)F(y+\rho^1,x+\rho^2))dxdy.
\end{eqnarray}
\end{lem}
\noindent
{\bf Proof of Lemma \ref{l:3.4}} From Lemma \ref{l:3.3} it follows that 
$\langle\Lambda(x_j+\rho^j,x_k+\rho^k;T), \UNO_{[0,t]}\rangle$ are well defined, and 
$\langle\Lambda(\cdot +\rho^j,\cdot +\rho^k;T), \UNO_{[0,t]}\rangle$ belongs to $L^2(\mathbb R^2\times \Omega,\lambda\otimes\lambda\otimes P)$,
hence by Lemma \ref{l:6.2} the process
 $\xi^T$ is well defined and the series in 
(\ref{eq:1.1}) converges in $L^2(\Omega)$. Moreover, using the fact that
$\Lambda(x+\rho^j,y+\rho^k,T)=\Lambda(y+\rho^k,x+\rho^j;T)$
(see Corollary 3.4 in \cite{BT1}), we have
\begin{equation}
\label{eq:6.11}
E\left(\sum_{j, k\atop j\neq k}\sigma_j\sigma_k\langle\Lambda(x_j+\rho^j,x_k+\rho^k;T),\psi\rangle\right)^2
=2\int_{\mathbb R^2}E\langle\Lambda(x+\rho^1,y+\rho^2;T),\psi\rangle^2dxdy,
\end{equation}
for $\psi$ either of the form (\ref{eq:3.6}), or $\psi=\psi_1+\varphi, \psi_1$
of the form (\ref{eq:3.6}), $\varphi\in{\cal S}$. Hence for $t_2>t_1\geq 0$, by 
(\ref{eq:1.1}), (\ref{eq:6.11}) and (\ref{eq:6.9}),
$$E\left(\xi^T_{t_2}-\xi^T_{t_1}\right)^2=\frac{1}{2\pi^2T^2}\int_{[0,T]^4}
\int_{\mathbb R^2}\left|\widehat{\UNO_{(t_1,t_2]}}(x+y)\right|^2e^{-|s-u||x|^\alpha}e^{-|r-v||y|^\alpha}dxdydsdrdudv.$$
Using $1>\alpha>\frac{1}{2}$ and
\begin{equation}
\label{eq:6.12}
\frac{1}{T}\int^T_0\int^T_0e^{-|s-r||x|^\alpha}dsdr\leq \frac{2}{|x|^\alpha}
\end{equation}
we obtain
\begin{eqnarray}
E(\xi_{t_2}^T-\xi^T_{t_1})^2&\leq&C\int_{\mathbb R^2}\frac{|e^{i(t_2-t_1)(x+y)}-1|^2}{|x+y|^2}|x|^{-\alpha}|y|^{-\alpha}dxdy\nonumber\\
\label{eq:6.13}
&\leq&C_1|t_2-t_1|^{2\alpha}.
\end{eqnarray}
Hence $\xi^T$ has a continuous modification. $\hfill\Box$
\vglue.5cm

Before we pass to the proof of Theorem \ref{t:3.5} we observe that for any ${\cal S}'$-random variable $Z$ (not necessarily Gaussian) such that ${\cal S}\ni\varphi
\mapsto E\langle Z,\varphi\rangle^2$ is continuous, the Wick product $:Z\otimes Z:$ is well-defined by an extension of (\ref{eq:3.1}). Moreover, we have the following lemma.

\begin{lem}
\label{l:6.3}
Let $(Z_T)_{T\geq 1}$ be a family of ${\cal S}'$-random variables such that
$$\sup_{T\geq 1}E\langle Z_T,\varphi\rangle^2\leq p^2(\varphi),\quad\varphi\in S,$$
form some continuous Hilbertian seminorm $p$ on ${\cal S}$. Assume that $Z_T\Rightarrow Z$ and $E\langle Z_T,\varphi\rangle^2\to E\langle Z,\varphi\rangle^2,\varphi\in{\cal S}$, as $T\to\infty$. Then $:Z_T\otimes Z_T:, :Z\otimes Z:$ are well defined and
$:Z_T\otimes Z_T:\Rightarrow \\
:Z\otimes Z:$  in ${\cal S}'(\mathbb R^2)$ as $T\to\infty$.
\end{lem}

This lemma follows by a standard argument using properties of ${\cal S}$ \cite{Tr}, so we skip the proof. 

Lemma \ref{l:6.3} together with Proposition \ref{p:2.1}(a) imply the following corollary.

\begin{cor}
\label{c:6.4}
Let $X_T,X$ be as in Proposition \ref{p:2.1}.
Then 
\begin{equation}
\label{eq:6.14}
:X_T\otimes X_T:\Rightarrow :X\otimes X:\, \hbox{\it in}\ {\cal S}'(\mathbb R^2)\;\;{\it as}
\;\; T\to\infty.
\end{equation}
\end{cor}

Indeed, it suffices to observe that by (\ref{eq:2.14}) and the Poisson initial condition we have
\begin{eqnarray*}
E\langle X_T,\varphi\rangle^2&=&\frac{1}{T}\int_{\mathbb R}\int^T_0\int^T_0E\varphi(x+\rho_s)\varphi(x+\rho_u)dudsdx\\
&=&\frac{1}{2\pi}\frac{1}{T}\int_{\mathbb R}\int^T_0\int^T_0|\widehat{\varphi}(x)|^2e^{-|s-u||x|^\alpha}dudsdx,
\end{eqnarray*}
so the assumptions of Lemma \ref{l:6.3} are satisfied (we use (\ref{eq:6.12})).
\vglue.5cm
\noindent
{\bf Proof of Theorem \ref{t:3.5}} For $\psi$ of the form (\ref{eq:3.6}) we denote by $\xi^T_\psi$ the random variable defined by (\ref{eq:1.1}) with $\UNO_{[0,t]}$ replaced by $\psi$.

To prove the theorem it suffices to show that
\begin{equation}
\label{eq:6.15}
\lim_{T\to\infty}|E e^{i\xi^T_\psi}-Ee^{i\langle Y,\psi\rangle}|=0
\end{equation}
for any $\psi$ of the form (\ref{eq:3.6}). Indeed, from (\ref{eq:6.15}) and Theorem \ref{p:3.2} we infer convergence of finite-dimensional distributions, and from 
(\ref{eq:6.13}) we obtain tightness in $C([0,\tau])$ for each $\tau>0$ (see \cite{Billingsley}, Thm. 12.3; note that the constant $C_1$ in (\ref{eq:6.13}) does not depend on $T$).

Fix any $f\in{\cal F}$ and denote $\psi_\delta=\psi * f_\delta,\delta>0$, where $f_\delta$ is given by (\ref{eq:2.15}). Let
\begin{equation}
\label{eq:6.16}
\xi^{T,\varepsilon,f}_\varphi=\frac{1}{T}\sum_{j,k\atop j\neq k}\sigma_j\sigma_k\langle\Lambda^f_\varepsilon(x_j+\rho^j,x_k+\rho^k;T),\varphi\rangle,\quad \varepsilon >0,\varphi\in{\cal S},
\end{equation}
which is well defined by Lemma \ref{l:6.2}.

Using the estimate $|Ee^{i\eta_1}-Ee^{i\eta_2}|\leq \frac{1}{2} E|\eta_1-\eta_2|^2,$ valid for centered random variables $\eta_1,\eta_2$, it is easy to see that 
(\ref{eq:6.15}) will be proved if we show
\begin{eqnarray}
\label{eq:6.17}
&&\lim_{\delta\to 0}\sup_{T\geq 1}E|\xi^T_\psi-\xi^T_{\psi_\delta}|^2=0,\\
\label{eq:6.18}
&&\lim_{\varepsilon\to 0}\sup_{T\geq 1}\sup_{0<\delta\leq 1}E|\xi^T_{\psi_\delta}-\xi^{T,\varepsilon,f}_{\psi_\delta}|^2=0,\\
\label{eq:6.19}
&&\lim_{T\to\infty}\sup_{0<\delta\leq 1}E|\langle : X_T\otimes X_T:,\Phi^f_{\varepsilon,\psi_\delta}\rangle-\xi^{T,\varepsilon,f}_{\psi_\delta}|=0,\quad \varepsilon >0,\\
\label{eq:6.20}
&&\langle :X_T\otimes X_T:,\Phi^f_{\varepsilon,\psi_\delta}\rangle\Rightarrow\langle :X\otimes X:,\Phi^f_{\varepsilon,\psi_\delta}\rangle\quad{\rm as}\quad T\to\infty,\quad \varepsilon>0,\delta>0,\\
\label{eq:6.21}
&&\lim_{\varepsilon\to 0}\sup_{0<\delta\leq 1}E|\langle :X\otimes X:,\Phi^f_{\varepsilon,\psi_\delta}\rangle-\langle Y,\psi_\delta\rangle|^2=0,\\
\label{eq:6.22}
&&\lim_{\delta\to 0} E|\langle Y,\psi_\delta\rangle-\langle Y,\psi\rangle|^2=0.
\end{eqnarray}

Using (\ref{eq:6.10}), (\ref{eq:6.9}), (\ref{eq:6.7}) and then (\ref{eq:6.12}), we have
\begin{eqnarray}
E|\xi^T_\psi-\xi^T_{\psi_\delta}|^2&=&\frac{2}{T^2}\int_{[0,T]^4}
\int_{\mathbb R^2}|\widehat{\psi}(x+y)|^2|1-\widehat{f}(\delta(x+y))|^2
\,\,e^{-|s-u||x|^\alpha}e^{-|r-v||y|^\alpha}dxdydsdudrdv\nonumber \\
\label{eq:6.23}
&\leq&8\int_{\mathbb R^2}|\widehat{\psi}(x+y)|^2|1-\widehat{f}(\delta(x+y))|^2|x|^{-\alpha}|y|^{-\alpha}dxdy.
\end{eqnarray}
Hence (\ref{eq:6.17}) follows by (\ref{eq:2.16}) and since $\frac{1}{2}<\alpha<1.$

Next, we apply Lemma \ref{l:6.2} to
$$F(x+\rho^1,y+\rho^2)=\langle\Lambda(x+\rho^1,y+\rho^2;T)-\Lambda^f_\varepsilon(x+\rho^1,y+\rho^2;T),\psi_\delta\rangle,$$
and by (\ref{eq:6.10}) we obtain
\begin{eqnarray*}
\lefteqn{E|\xi^T_{\psi_\delta}-\xi^{T,\varepsilon,f}_{\psi_\delta}|^2\leq \frac{2}{T^2}\int_{\mathbb R^2}E(\langle(\Lambda-\Lambda^f_\varepsilon)(x+\rho^1,y+\rho^2;T),\psi_\delta\rangle)^2dxdy}\\
&\leq&\frac{C}{T^2}\int_{[0,T]^4}\int_{\mathbb R^2}|\widehat{\psi}(x+y)|^2|\widehat{f}(\varepsilon y)-1|^2e^{-|s-u||x|^\alpha}e^{-|r-v||y|^\alpha}dxdydsdrdudv,
\end{eqnarray*}
by (\ref{eq:6.8}) and (\ref{eq:2.16}). By (\ref{eq:6.12}), (\ref{eq:6.3}) and (\ref{eq:2.16}) we obtain (\ref{eq:6.18}).

To prove (\ref{eq:6.19}), observe that random variables almost identical to $\xi^{T,\varepsilon,f}_{\psi_\delta}$ and \\
$\langle :X_T\otimes X_T:,\Phi^f_{\varepsilon,\psi_\delta}\rangle$ have already appeared in \cite{BT1} with different notation (and different scaling). Hence, by (4.1), (8.8), (8.14) and the two subsequent formulas in \cite{BT1}, we have
\begin{eqnarray}
A(T,\varepsilon,\delta)&:=&E|\langle :X_T\otimes X_T:,\Phi^f_{\varepsilon,\psi_\delta}\rangle-\xi^{T,\varepsilon,f}_{\psi_\delta}|^2\nonumber\\
\label{eq:6.24}
&=&\frac{1}{T^2}\int_{[0,T]^4}\int_{\mathbb R}E\Phi^f_{\varepsilon,\psi_\delta}(x+\rho_s,x+\rho_u)\Phi^f_{\varepsilon,\psi_\delta}(x+\rho_r,x+\rho_v)dxdvdrduds.
\end{eqnarray}
It is easy to see  that since the support of $\psi_\delta$ is contained in a compact set which is independent of $\delta$, we have (see (\ref{eq:3.2}) and (\ref{eq:3.6}))
\begin{equation}
\label{eq:6.25}
|\Phi^f_{\varepsilon,\psi_\delta}(x,y)|\leq C(\varepsilon,f)\phi(x)\phi(y),
\end{equation}
where
$$\phi(x)=\frac{1}{1+|x|^2}.$$

Let $({\cal T}_t)_t$ denote the $\alpha$-stable semigroup and $G$  its potential, $G\varphi=\int^\infty_0{\cal T}_t\varphi dt$.

From (\ref{eq:6.24}), (\ref{eq:6.25}), and the Markov property we get
\begin{eqnarray*}
\lefteqn{A(T,\varepsilon,\delta)}\\
&\leq&4! C^2(\varepsilon,f)\frac{1}{T^2}\int_{\mathbb R}\int^T_0\int^T_s\int^T_u\int^T_r{\cal T}_s(\phi({\cal T}_{u-s}\phi({\cal T}_{r-u}(\phi{\cal T}_{v-r}\phi))))(x) dxdvdrduds\\
&\leq&C_1(\varepsilon,f)\frac{1}{T^2}\int^T_0\int_{\mathbb R}\phi(x)G(\phi G(\phi G\phi))(x)dxds\\
&\leq&C_2(\varepsilon,f)\frac{1}{T}.
\end{eqnarray*}
In the last estimate we used the fact that $G\phi$ is bounded. Hence (\ref{eq:6.19}) follows.

(\ref{eq:6.20}) follows from (\ref{eq:6.14}) since $\Phi^f_{\varepsilon,\varphi_\delta}\in{\cal S}$.

In the proof of (\ref{eq:6.21}) we use (\ref{eq:3.3}) and Lemma \ref{l:3.1}, which yield
\begin{eqnarray*}
\lefteqn{E|\langle :X\otimes X:,\Phi^f_{\varepsilon,\psi_\delta}\rangle-\langle Y,\psi_\delta\rangle|^2}\\
&=&\frac{1}{\pi^2}\int_{\mathbb R^2}|\widehat{\psi}_\delta(x+y)|^2|\widehat{f}(\varepsilon y)-1|^2|x|^{-\alpha}|y|^{-\alpha}dxdy\\
&&+\frac{1}{\pi^2}\int_{\mathbb R^2}|\widehat{\psi}_\delta(x+y)|^2
(\widehat{f}(\varepsilon y)-1)(\overline{\widehat{f}(\varepsilon x)-1})
|x|^{-\alpha}|y|^{-\alpha}dxdy\\
&\leq&\frac{2}{\pi^2}\int_{\mathbb R^2}|\widehat{\psi}(x+y)|^2|\widehat{f}(\varepsilon y)-1|^2|x|^{-\alpha}|y|^{-\alpha}dxdy.
\end{eqnarray*}
Hence (\ref{eq:6.21}) follows by (\ref{eq:6.3}) and (\ref{eq:2.16}).

Finally, from (\ref{eq:3.3}) and Lemma \ref{l:3.1} it is easy to see that
$$E|\langle Y,\psi_\delta\rangle-\langle Y,\psi\rangle|^2=\frac{2}{\pi^2}\int_{\mathbb R}|\widehat{\psi}(x+y)|^2|\widehat{f}(\delta(x+y))-1|^2|x|^{-\alpha}|y|^{-\alpha}dxdy,$$
which tends to $0$ as $\delta\to 0$ (cf (\ref{eq:6.23})). $\hfill\Box$
\vglue.5cm
\noindent
{\bf Proof of Theorem \ref{t:4.1}} Fix $0\leq s<t$. By (\ref{eq:2.1}) and (\ref{eq:2.2})
we know that
\begin{equation}
\label{eq:6.26}
Ee^{i[z_1(\xi_t-\xi_s)+z_2(\xi_{t+\tau}-\xi_{s-\tau})]}
={\rm exp}\left\{\frac{1}{2}\sum^\infty_{k=2}\frac{(2i\sigma)^k}{k}\widetilde{R}_k(z_1,z_2,s,t,\tau)\right\},
\end{equation}
where $\widetilde{R}_k(z_1,z_2,s,t,T)$ is the right hand side of (\ref{eq:2.2}) with
\begin{equation}
\label{eq:6.27}
\psi(x)=z_1\UNO_{(s,t]}+z_2\UNO_{(s+\tau,t+\tau]},
\end{equation}
and formula (\ref{eq:6.26}) holds for $z_1,z_2\in{\mathbb R}$ and $\tau>0$ such that the series in (\ref{eq:6.26}) converges.

To continue with the proof we need the following lemma and corollary.

\begin{lem}
\label{l:6.5}
There exists $C(s,t)>0$ independent of $\tau$ such that
\begin{equation}
\label{eq:6.28}
|\widetilde{R}_k(z_1,z_2,s,t,\tau)|\leq (|z_1|+|z_2|)^k(C(s,t))^k,\quad k=2,3,\ldots .
\end{equation}
\end{lem}

An immediate consequence of this lemma is 

\begin{cor}
\label{c:6.6}
For all $0<s<t$ there exists a neighborhood $U(s,t)$ of $0$ in $\mathbb R^2$ such that (\ref{eq:6.26}) holds for $(z_1,z_2)\in U(s,t)$ and all $\tau>0$.
\end{cor}
\noindent
{\bf Proof of Lemma \ref{l:6.5}} Consider an integral
\begin{equation}
\label{eq:6.29}
I_k(\varepsilon_1,\ldots,\varepsilon_k)=\int^{t+\varepsilon_1\tau}_{s+\varepsilon_1\tau}\ldots\int^{t+\varepsilon_k\tau}_{s+\varepsilon_k\tau}|x_1-x_2|^{H-1}\ldots |x_{k-1}-x_k|^{H-1}|x_k-x_1|^{H-1}dx_1\ldots dx_k,
\end{equation}
where $(\varepsilon_1,\ldots,\varepsilon_k)\in\{0,1\}^k$. For each $j$ we substitute $x'_j=x_j-\varepsilon_j\tau$ and use
$$|y-x+\tau|^{H-1}\leq |y-x|^{H-1}\quad{\rm for}\quad x,y\in[s,t],\tau >\tau_0(s,t).$$
Hence, similarly as in formula (13) in \cite{T4},
\begin{equation}
\label{eq:6.30}
I_k(\varepsilon_1,\ldots,\varepsilon_k)\leq (C(s,t))^k.
\end{equation}
The same inequality holds for $\tau\leq \tau_0(s,t)$.

By (\ref{eq:6.27}), (\ref{eq:2.2}), 
\begin{equation}
\label{eq:6.31}
\widetilde{R}_k(z_1,z_2,s,t,\tau)=z^k_1I_k(0,\ldots,0)+\widetilde{\widetilde{R}}_k(z_1,z_2,s,t,\tau)+z^k_2I_k(0,\ldots,0),
\end{equation}

\noindent
where
\begin{equation}
\label{eq:6.32}
\widetilde{\widetilde{R}}_k(z_1,z_2,s,t,\tau)=\sum^{k-1}_{j=1}z^{k-j}_1z^j_2
\sum_{(\varepsilon_1,\ldots,\varepsilon_k)\in\{0,1\}^k\atop \varepsilon_1+\ldots+\varepsilon_k=j}I_k(\varepsilon_1,\ldots,\varepsilon_k).
\end{equation}
Hence (\ref{eq:6.30})
  implies 
\begin{equation}
\label{eq:6.33}
|\widetilde{R}_k(z_1,z_2,s,t,\tau)|\leq\sum^k_{j=1}
\left({k\atop j}\right)|z_1|^{k-j}|z_2|^j(C(s,t))^k=(|z_1|+|z_2|)^k(C(s,t))^k,
\end{equation}
which proves the lemma. $\hfill\Box$
\vglue .25cm
We return to the proof of the theorem.

Let $U(s,t)$ be as in Corollary \ref{c:6.6}. It is now clear that for $(z_1,z_2)$ in $U(s,t)$ the function $D^\xi_\tau(z_1,z_2,s,t)$ given by (\ref{eq:4.3}) is well defined for all $\tau>0$.

Moreover, by (\ref{eq:4.3}), (\ref{eq:4.2}), (\ref{eq:6.26}) and (\ref{eq:6.31}),
\begin{equation}
\label{eq:6.34}
D^\xi_\tau(z_1,z_2,s,t)=\left|\frac{1}{2}\sum^\infty_{k=2}\frac{(2i\sigma)^k}{k}\widetilde{\widetilde{R}}_k(z_1,z_2,s,t,\tau)\right|
\end{equation}

It is not difficult to see that due to the circular form of the integrand in each $I_k$ in (\ref{eq:6.32}), after the substitution $x'_j=x_j-\varepsilon_j\tau$ (see (\ref{eq:6.29})) there appear at least two factors of the form $|y-x+\tau|^{H-1}, y,x\in[s,t]$, and for some $(\varepsilon_1,\ldots,\varepsilon_k)$ there are exactly two such factors. Hence, for each $\gamma<2-2H$,
$$\lim_{\tau\to\infty}\tau ^\gamma\widetilde{\widetilde{R}}_k(z_1,z_2,s,t,\tau)=0,$$
and for some $z_1,z_2$
$$\lim_{\tau\to\infty}\tau^{2-2H}\widetilde{\widetilde{R}}_u(z_1,z_2,s,t,\tau)\neq 0.$$
Hence the theorem follows by Lemma \ref{l:6.5}. $\hfill\Box$
%

\end{document}